\input amstex.tex
\baselineskip 22pt
\documentstyle{amsppt}
\topmatter
\title
Some remarks on the universal cover of an open K3 surface
\endtitle
\author
Fabrizio Catanese, JongHae Keum and Keiji Oguiso
\endauthor 
\address
Mathematisches Institut der Georg-August Universit\"at
G\"ottingen, 
\flushpar 
Bunsen-strasse 3-5 D-37073 G\"ottingen, Germany,
\flushpar
Email : catanese\@uni-math.gwdg.de
\endaddress
\address
Korea Institute for Advanced Study, 207-43 Cheongryangri-dong,
Dongdaemun-gu, Seoul 130-012, Korea, Email : jhkeum\@kias.re.kr
\endaddress
\address
Department of Mathematical Sciences, University of Tokyo, Komaba,
Meguro, Tokyo, Japan, Email : oguiso\@ms.u-tokyo.ac.jp
\endaddress

\subjclass
14J28
\endsubjclass

\abstract We shall give, in an optimal form, a sufficient numerical  condition
for the finiteness of the fundamental group of the smooth locus of a normal
K3 surface. We shall moreover prove that, if the normal K3 surface is elliptic
and the above fundamental group is not finite, then there is a finite covering
which is a complex torus.
\endabstract

\leftheadtext{F. Catanese, J. Keum and K. Oguiso}
\rightheadtext{Universal cover of an open K3 surface}
\endtopmatter

\document
\head
{Introduction}
\endhead

Throughout this note, we work in the category of separated complex analytic
spaces  endowed with the  classical topology.

\par
\vskip 4pt

The first aim of this short note is to add some evidence to the following
interesting Conjecture posed by De-Qi Zhang (See [KT], [SZ], [KZ] and [Ca1,2]
for related work):

\proclaim{Conjecture} The universal cover of the smooth locus of
a normal K3 surface is  a big open set of either a normal K3 surface or  of
$\bold C^{2}$. Furthermore, in the latter case, the universal cover factors
through a finite \'etale cover by a big open set of a torus.
\endproclaim

To explain our notation, here and hereafter, a {\it normal K3 surface} means a
normal surface whose minimal resolution is a K3 surface and a {\it torus} is a
2-dimensional complex torus. By a {\it big open set} we mean the
complement of a discrete subset.

Note  that, by a result of Siu [Si], a K3 surface is always
K\"ahler.

Since the canonical divisor of a smooth K3 surface is trivial, it follows that
the singularities of a normal K3 surface can only be the so called Du Val
singularities ( also called Rational Double Points): these are also the
Kleinian singularities obtained as quotients $\bold C^{2}/ G$ with $G
\subset SL(2,\bold C$).

The property that the singularities $x_i$ are of the form $\bold C^{2}/ G_i$
allows to define the orbifold Euler number of the normal K3 surface.
For each singular point $x_i$ we take a neighbourhood $U_i$ of $x_i$ which is
the quotient of a ball in $\bold C^{2}$ and decree that its orbifold Euler
number equals $1/|G_i|$: using a Mayer Vietoris sequence we can then extend
the definition to the whole of the normal K3 surface.
One has that the orbifold Euler number is non negative, and indeed
 R. Kobayashi and  A. Todorov [KT] showed that the second case in
Zhang's conjecture happens if and only if the orbifold  Euler number of the
normal K3 surface is zero.

\par
\vskip 4pt Let $\overline{S}$ be a normal K3 surface, let $S^{0}
:= \overline{S} - \text{Sing}\, \overline{S}$ and let $\nu : S
\rightarrow \overline{S}$ the minimal resolution.

Our first result is

\proclaim{Theorem A} If the normal K3 surface $\overline{S}$
admits an elliptic fibration then either $\pi_{1}(S^{0})$ is
finite or there is a finite covering of $\overline{S}$, ramified
only on a finite set, which yields a complex torus.
\endproclaim
\par
\vskip 4pt

Our second aim in this note is to establish a sharp sufficient condition for
the validity of the first alternative in Zhang's conjecture.

We set $E$ to be the reduced exceptional divisor $E :=
\nu^{-1}(\text{Sing}\, \overline{S})$ and decompose $E$ into
irreducible components $E := \sum_{i = 1}^{r} E_{i}$. An important invariant
is the number $r := r(\overline{S})$  of
irreducible components of the exceptional divisor. Clearly, $S^{0} = S -
E$.
\par
\vskip 4pt

Our main observation is as follows:

\proclaim{Theorem B} If $r = r(\overline{S}) \leq 15$, then
$\pi_{1}(S^{0})$ is finite. In particular, if $r \leq 15$, then
the universal cover of $S^{0}$ is a big open set of a normal K3
surface. \endproclaim
\par
\vskip 4pt

This easy remark however gives us the best possible uniform
bound on $r$ in order that $\pi_{1}(S^{0})$  be finite, in view of
the following facts [KT], [KZ]:

\roster
\item A normal Kummer surface $A/-1$ satisfies
$r = 16$ and $\vert \pi_{1}((A/-1)^{0}) \vert = \infty$.
\item There is a normal K3 surface $\overline{S}$ such that
$r(\overline{S}) \le 15$ but $\pi_{1}(S^{0}) \not= \{1\}$.
\endroster

In our actual proof, the numerical condition $r \leq 15$ will be used
twice: in Lemmas 1 and 2. These two Lemmas allow us to  establish the
existence of an elliptic pencil, whence to reduce the proof of theorem B to
theorem A.
 For the last statement of Theorem B, we recall that the category of Du Val
singularities is closed under the operation of taking the normalization of a
finite covering which is unramified outside a finite set.

Therefore, by the
classification of smooth compact K\"ahler surfaces
due to Castelnuovo, Enriques, Kodaira, the normalization of a finite cover of
a normal K3 surface is either a normal K3 surface or a 2-dimensional torus if
the covering
is unramified outside a finite set . This fact will be frequently used in
our proof, too.

The proof of theorem A is on one side  based on an exact sequence for
open fibred surfaces, and which
relies on the notion of orbifold fundamental group (cf. e.g. [Ca1]).

On the other hand, it is based on the existence of finite branched covers of
given branching type provided by finite index subgroups of the orbifold
fundamental group.


\par
\vskip 4pt Of course,  the most interesting part of the
Conjecture concerns the case where the fundamental group is infinite.
As such, the question seems to belong more to the transcendental theory and we
hope to return on the question using the current technologies on uniformization
problems.

In the last paragraph
of this note (Remark 4) we will
just comment on a more or less obvious reduction process.

\par
\vskip 4pt
\head Acknowledgement \endhead

An initial idea of this work was found during the first and the third authors'
stay at KIAS in Seoul under the financial support by the Institute.
They would like to express their gratitude to Professor Jun-Muk Hwang
and KIAS for making their stay enjoyable. The third author would like
to express his thanks to Professors Yujiro Kawamata and Eckart Viehweg for
their warm encouragement during his stay at the Institute.
\par
\vskip 4pt

\head Proofs \endhead

\proclaim{Lemma 1} If $r = r(\overline{S}) \leq 15$ then $\overline{S}$ does
not admit any finite covering by a complex torus which is unramified outside a
finite set.
 \endproclaim

\demo{Proof}

Assuming  the contrary, we shall show that $r \geq
16$. Let $Q_{k}$ ($1 \leq k \leq m)$ be the singular points of
$\overline{S}$ and let $n_{k}$  be the number of the irreducible
components of $\nu^{-1}(Q_{k})$. Note that the contribution of
$Q_{k}$ to the orbifold Euler number of $S$ is $(n_{k} + 1) -
1/\delta_{k}$, where $\delta_{k}$ is the order of the local
fundamental group around $Q_{k}$. Since the orbifold Euler number of
$\overline{S}$ equals $ 0$ (the obvious
direction of [KT]) one has

$$ e(S) = 24 = \sum_{k = 1}^{m} (n_{k} + 1 - \frac{1}{\delta_{k}}).$$
Note that $\delta_{k} = 1/2$ if $n_{k} = 1$. Then, one has also
$$n_{k} + 1  - \frac{1}{\delta_{k}} \leq \frac{3}{2}n_{k}.$$

Now, substituting the second inequality into the first equality,
one obtains
$$24 \leq \frac{3}{2} \sum_{k = 1}^{m} n_{k} = \frac{3}{2}r.$$
This implies $r \geq 16$. \qed \enddemo

\proclaim{Lemma 2} If $r = r(\overline{S}) \leq 15$ then there exist an
elliptic K3 surface
$f' : S' \rightarrow \bold P^{1}$ and an effective divisor $E'$ supported
in fibers of $f'$ such that $S^{0}$ is diffeomorphic to $S' - E'$ and such that
$E$ and $E'$ are of the same Dynkin type.
\endproclaim

\demo{Proof}

We argue by  descending induction on $\rho(S)$.
If $\rho (S) = 20$, then $S$ is algebraic and the
orthogonal complement of
$\bold Z\langle E_{i} \vert 1 \leq i \leq r \rangle$ in $\text{Pic}(S)$
is a hyperbolic lattice of rank $20 - r(S) \geq 20 - 15 = 5$.

Then, by Meyer's theorem (cf. e.g. [Se], Chapter IV, 3, Cor. 2),
there exists $X \in \text{Pic}(S)- \{0\}$ such that $(X.E_{i}) =
0$ for all $E_{i}$ and such that $(X^{2}) = 0$. ({\it However, in
general this} $X$ {\it does not lie in the nef cone of} $S$.)

Let $u : \Cal U \rightarrow \Cal K$ be the Kuranishi
family of $S$ and choose a trivialization $R^{2}u_{*}\bold Z_{\Cal U} \simeq
\Lambda$ over $\Cal K$. Here $\Lambda$ is an even unimodular lattice
of index $(3, 19)$ called the K3 lattice. Then one can define the period
map $p :
\Cal K \rightarrow \Cal P$, where $\Cal P$ is the period domain.
This $p$ is a local isomorphism by the local Torelli Theorem and
allows us to identify $\Cal K$ with a small open set of the
period domain $\Cal P$ (denoted again by $\Cal P$).

Let's denote by $e_{i}$, $x$
the elements in $\Lambda$ corresponding to
$E_{i}$ and $X$. Let us take
the sublocus $\Cal B \subset \Cal K$ defined by the equations

$$\{ \omega \in {\Cal P}| (e_{1}. \omega) = (e_{2}. \omega) = \cdots = (e_{r}.
\omega) = (x. \omega) = 0\},$$

and consider the induced family $\Cal T \rightarrow \Cal B$.
Here $\Cal B$ is of dimension
$20 - 1 - r > 0$. Note that $h^{0}(\Cal O_{S}(E_{i})) = 1$ and
$h^{1}(\Cal O_{S}(E_{i})) = 0$.

Then, by  construction and by the base change Theorem, the smooth
rational curves $E_{i}$ lift uniquely to  effective divisors
$\Cal E_{i}$ flat over $\Cal B$ in such a way that the fibers
$\Cal E_{i, b}$ are smooth rational curves (of the same Dynkin
type as $E_{i}$) and $S - \cup_{i = 1}^{r} E_{i}$ is
diffeomorphic to $\Cal T_{b} - \cup_{i = 1}^{r} \Cal E_{i, b}$ for
all $b \in \Cal B$ if $\Cal B$ is chosen to be sufficiently
small. Furthermore, for generic $b \in \Cal B$, the Picard group
of $\Cal T_{b}$ is isomorphic to the primitive closure of $\bold
Z\langle x, e_{i} \vert 1 \leq i \leq r \rangle$ in $\Lambda$
(See for instance [Og]) and is a semi-negative definite lattice.
Thus $\Cal T_{b}$ is of algebraic dimension one. Now, the
algebraic reduction map of this $\Cal T_{b}$ gives an elliptic
fibration with the required properties.
\par
\vskip 4pt

Assume that $\rho := \rho (S) < 20$. As before, we let
$u : \Cal U \rightarrow \Cal K$ the Kuranishi family of $S$
and choose a trivialization $R^{2}u_{*}\bold Z_{\Cal U} \simeq
\Lambda$ over $\Cal K$ and identify $\Cal K$ with a small open set of the
period domain $\Cal P$. By $x_{i}$
($1 \leq i \leq \rho$) we denote the elements in $\Lambda$
corresponding to some integral basis of $\text{Pic}(S)$. Let us
take the sublocus $\Cal A \subset \Cal K = \Cal P$ defined by the
equations

$$\{ \omega \in {\Cal P}| (x_{1}. \omega) = (x_{2}. \omega) = \cdots =
(x_{\rho}. \omega) = 0 \},$$
and consider the induced family $\pi :
\Cal S \rightarrow \Cal A$. Here $\Cal A$ is of dimension $20 -
\rho > 0$.

Then, for the same reason as before, the smooth rational curves
$E_{i}$ lift
uniquely to  effective divisors $\Cal C_{i}$ flat over $\Cal
A$ in such a way that the fibers $\Cal C_{i, a}$ are smooth
rational curves (of the same Dynkin type as $E_{i}$) and $S -
\cup_{i = 1}^{r} E_{i}$ is diffeomorphic to $\Cal S_{a} - \cup_{i
= 1}^{r} \Cal C_{i, a}$ for all $a \in \Cal A$ if $\Cal A$ is
chosen to be sufficiently small. Therefore the orbifold Euler
numbers are also the same for $\overline{S}$ and $\overline{\Cal
S}_{a}$. Here $\overline{\Cal S}_{a}$ is the normal K3 surface
obtained by the contraction of $\cup_{i = 1}^{r} \Cal C_{i, a}$.
Furthermore, by  construction, $\pi$ is not isotrivial and the fibers
$\Cal S_{a}$
satisfy $\rho(\Cal S_{a}) \geq \rho$. Then by [Og] there is $a
\in \Cal A$ such that $\rho (\Cal S_{a}) > \rho$. Now,
we are done by  descending
induction on $\rho(S)$. \qed \enddemo
\par
\vskip 4pt

By Lemmas 1 and 2, it is now clear that Theorem A implies Theorem
B. Before  we proceed to the proof of theorem A we state a very
easy and quite general lemma

\proclaim{Lemma 3}
Let  $f : S \rightarrow C$ be a fibration of a complete surface onto a
complete curve   with general fiber $F$, and let $S^{0}$ be the Zariski open
 set of $S$ which is defined as the complement of a finite number of divisors
 $Z_{j}$ each contained in a fibre $F_{j} $, but with support properly
contained
in $ supp(F_{j}) $.

Denote by $W_{j}$ the maximal divisor $ \leq F_{j}$ which has no common
component with $Z_{j}$, and by $m_{j}$ the greatest common divisor of the
multiplicities of the components of $W_{j}$.

Set $P_{j} = f(F_{j}) $, assume that the points $\{P_{1},
\cdots , P_{k}\}$ are all distinct, denote by $B :=
C - \{P_{1}, \cdots , P_{k} \}$ and define the orbifold fundamental group $
\pi_1^{orb}(f^{0}) $ of the open fibration $f^{0} : S^{0} \rightarrow C$ as the
quotient of
$ \pi_1(B) $ by the normal subgroup generated by the elements
$ \gamma_j^{m_{j}}$, where $ \gamma_j$ is a simple loop going around the point
$P_j$. Then the fundamental group $ \pi_1(S^{0}) $ fits into an exact sequence

$$\pi_{1}(F)  \rightarrow \pi_{1}(S^{0}) \rightarrow
\pi_1^{orb}(f^{0}) \rightarrow 1.$$
\endproclaim

\demo{Proof}
Let $S'$ be the complement of the given fibres, i.e., $S' = S -
\{F_{1}, \cdots , F_{k}\}$, so that $ f' : S' \rightarrow B$ is a fibre bundle
with fibre $F$.
Then we have the homotopy exact sequence

$$\pi_{1}(F)  \rightarrow \pi_{1}(S') \rightarrow
\pi_{1}(B) \rightarrow 1.$$

However, $S' = S^{0} - \{D_{1}, \cdots , D_{k}\}$ where $D_{j} :=
F_{j} \vert S^{0}$: therefore the kernel $K$ of the surjection
$\pi_{1}(S')  \rightarrow \pi_{1}(S^{0})$ is normally generated
by simple loops $\delta_{j,i}$ going around the components
$D_{j,i}$ of the $D_{j}$'s.

Each $\delta_{j,i}$ maps into $\pi_{1}(B)$ to the $m_{j,i}$-th
power of a conjugate of $\gamma_j$, where $m_{j,i}$ is the
multiplicity of $D_{j,i}$ in $D_{j}$. Whence the image $K'$ in
$\pi_{1}(B)$ of $K$ is normally generated by the elements
$\gamma_j^{m_{j,i}}$. Since $D_{j} = W_{j}\vert S^{0}$ we see
that $K'$  is normally generated by the elements $\gamma_j^{m_j}$,
thus we have the desired surjection $ \pi_{1}(S^{0}) \rightarrow
\pi_1^{orb}(f^{0}) \rightarrow 1$ with kernel generated by
$\pi_{1}(F)$. \qed
\enddemo
\proclaim{Theorem A}

Let $\overline{S}$ be a normal K3 surface such that
$\overline{S}$ admits an elliptic fibration

$\overline{f} : \overline{S} \rightarrow \bold P^{1}$. Then either
\roster
\item
 $\pi_{1}(S^{0})$ is finite or

\item $\overline{S}$  admits a finite covering  by a torus,
ramified only on a
finite set.
\endroster
\endproclaim

\demo{Proof}
 Note  by the way that all the assumptions  and assertions made in the  theorem
are stable  under replacement of $S^{0}$ by a finite  unramified covering.

\par
\vskip 4pt

Let $f : S \rightarrow \bold P^{1}$ be the elliptic fibration
induced by $\overline{f}$. We denote by $F$ a general fiber of
$f$ and set as in lemma 3 $f(E) := \{P_{1}, \cdots , P_{k}\}$, $B
:= \bold P^{1} - \{P_{1}, \cdots , P_{k} \}$, let $F_{j} :=
f^{-1}(P_{j})$ be the scheme theoretic fiber, and $D_{j} := F_{j}
\vert S^{0}$. Note that $\text{Supp}\, D_{j}$ is not empty for
each $j$ because the intersection product on the fibre is not
negative definite.

Lemma 3 gives us an exact sequence

$$\pi_{1}(F) \simeq \bold Z^{\oplus 2} \rightarrow \pi_{1}(S^{0}) \rightarrow
\pi_1^{orb}(f^{0}) \rightarrow 1.$$

We subdivide our analysis into two cases:

\roster
\item
 $\pi_1^{orb}(f^{0})$ is finite.
\item
 $\pi_1^{orb}(f^{0})$ is infinite.
\endroster

 In case (1), we take the unramified covering of $S^{0}$ associated
 with
 the epimorphism onto $\pi_1^{orb}(f^{0})$. We thus get another elliptic normal
K3 surface since the minimal resolution has the canonical divisor
trivial and the first Betti number $b_1 \leq 2$ ( by the
classification theorem, either $b_1=4 $ and we have a torus, or
$b_1=0 $ and we have a K3 surface), and we are reduced to the case
where $\pi_1^{orb}(f^{0})$ is trivial.

In this case it is immediate to see that, the fundamental group $\pi_{1}(F)
\simeq \bold Z^{\oplus 2}$ being abelian, the fundamental group
$\pi_{1}(S^{0})$
is the quotient of $\pi_{1}(F)
\simeq \bold Z^{\oplus 2}$ by the subgroup generated by the images of
$I_2 - T_j$, $I_2$ being the identity matrix and $T_j$ being the local
monodromy
matrix around the point $P_j$.

However, this quotient equals precisely the fundamental group of
the complete smooth surface $S$, which is a K3 surface. But then
$S$ is simply connected and we have proven our assertion.

In case (2), we use that $\pi_1^{orb}(f^{0})$ is the covering group of
a non compact simply connected Riemann surface $\Sigma$ branched over $ \bold
P^{1}$ with branching locus equal to $ \{P_{1}, \cdots , P_{k} \}$ and
branching
multiplicities $\{ m_{1}, \cdots , m_{k} \}$ (i.e., $\pi_1^{orb}(f^{0})$ acts
on $\Sigma$ with quotient $ \bold P^{1}$).

If the Riemann surface  $\Sigma$ were the disk then we would get a Fuchsian
group and there is a normal subgroup $\Gamma \subset \pi_1^{orb}(f^{0})$ of
finite index acting freely on $\Sigma$, whence a finite  Galois covering $C
\rightarrow \bold P^{1}$ where $C$ has genus at least $2$.

The epimorphism $\pi_{1}(S^{0}) \rightarrow \pi_1^{orb}(f^{0})/ \Gamma$ yields
an unramified covering of $S^{0}$ which compactifies to a smooth surface with
trivial canonical bundle. But such a surface admits no nontrivial holomorphic
map to a curve of genus $\geq 2$, thus we conclude that $\Sigma$ is the affine
line.

By considering again a normal subgroup $\Gamma \subset \pi_1^{orb}(f^{0})$ of
finite index acting freely on $\Sigma$, we get an elliptic curve $C$ and a
Galois covering $C \rightarrow \bold P^{1}$. Taking the normalization of
the fibre product we obtain again an unramified covering of $S^{0}$ which
compactifies to a surface
$X$ with trivial canonical bundle. Since we get a map  of $X$ to an elliptic
curve, $X$ is a torus (observe that $X$ is minimal), what proves our claims.
\qed
\enddemo

\remark{Remark 4} If $\pi_{1}(S^{0})$ is residually finite for
all $\overline{S}$ with $e_{\text{orb}}(\overline{S})
> 0$, then the Conjecture is affirmative.
\endremark

\demo{Proof} The case $e_{\text{orb}}(\overline{S}) = 0$ is
covered by [KT]. The remaining case follows from the assumption
together with the uniform boundedness of the order of finite
automorphism groups of K3 surfaces.  Here the uniform boundedness
is a consequence of the global Torelli Theorem together with the
Burnside property of finite groups of $GL_{n}(\bold Z)$. In the
algebraic case, one may use the result of Mukai [Mu] instead.
\qed \enddemo

\enddocument